\documentclass{ifacconf}

\usepackage[dvipdfmx]{graphicx}      
\usepackage{natbib}        
\usepackage{amsmath,amssymb,amsfonts}
\usepackage{mathrsfs}
\usepackage{algorithmic}
\usepackage[dvipdfmx]{graphicx}
\usepackage{color}

\allowdisplaybreaks

\begin{document}
\begin{frontmatter}

\title{Convergence Analysis of Natural Power Method and Its Applications to Control} 

\thanks[footnoteinfo]{This work was supported by JSPS KAKENHI Grant Numbers 23K26126, 21K12097, and 26K07550.}

\author[First]{Daiki Tsuzuki} 
\author[Second]{Kentaro Ohki} 

\address[First]{Graduate School of Informatics, Kyoto University, Yoshida-Honmachi, Sakyo-ku, Kyoto, Japan}
\address[Second]{Department of Applied Computer Engineering, Tokai University, Kitakaname 4-1-1, Hiratsuka, Kanagawa, Japan. (e-mail: ohki@tokai.ac.jp)}

\begin{abstract}                
This paper analyzes the discrete-time natural power method, demonstrating its convergence to the dominant $r$-dimensional subspace corresponding to the $r$ eigenvalues with the largest absolute values. This contrasts with the Oja flow, which targets eigenvalues with the largest real parts. We leverage this property to develop methods for model order reduction and low-rank controller synthesis for discrete-time LTI systems, proving preservation of key system properties. We also extend the low-rank control framework to slowly-varying LTV systems, showing its utility for tracking time-varying dominant subspaces.
\end{abstract}


\end{frontmatter}

\section{Introduction}
Online principal component algorithms have recently garnered significant attention. 
Oja's component flow (Oja flow), generated by a matrix differential equation introduced by \cite{oja1982simplified,oja1985stochastic}, is a prominent example. Its convergence properties have been investigated by \cite{chen1998global} for positive definite matrices, and by \cite{kuo2025asymptotic} and \cite{TsuzukiOhki2025global} for general square matrices. 
Discrete-time counterparts were subsequently proposed by \cite{oja1989neural} and \cite{hua1999new} primarily to extract dominant singular subspaces for positive definite matrices. 
Notably, Oja's discrete-time algorithm is the Euler scheme of the continuous-time Oja flow and extracts the same dominant subspace as its continuous counterpart.

\cite{hua1999new} proposed the following principal subspace tracking algorithm, also known as the \textit{natural power method (NPM)} or \textit{power iteration method}:
\begin{align}
    U[k+1] =& A U[k] (U[k]^{\top} A^{\top} A U[k]) ^{-1/2}, 
    \label{eq:natural_power_method}
\end{align}
where $U[0]$ is in the Stiefel manifold $\mathrm{St}(r,n)$ and $X^{1/2}$ denotes the square root of a positive semidefinite matrix $X$. 
Convergence analysis of this algorithm has primarily been investigated for cases where $A$ is a symmetric, positive-definite matrix. 
For applications involving the principal components of covariance matrices, the positive-semidefiniteness of $A$ is a natural requirement.
However, for other applications, such as extracting the left- or right-eigensubspaces of a non-symmetric matrix, positive semidefiniteness is not always appropriate.

This paper investigates the convergence properties of the discrete-time principal component algorithm \eqref{eq:natural_power_method} for general square matrices. 
We demonstrate that the NPM extracts the $r$-dominant subspace corresponding to the $r$ eigenvalues with the largest \textit{absolute values}. 
This result contrasts with the Euler approximation of the Oja flow, which extracts the $r$-dominant subspace corresponding to the $r$ eigenvalues with the largest \textit{real parts} (\cite{kuo2025asymptotic,TsuzukiOhki2025global}). 
Furthermore, we consider applications to control problems. 
Since the dominant subspace extracted by \eqref{eq:natural_power_method} corresponds to the dominant mode of discrete-time dynamical systems, the algorithm is naturally applicable to the model order reduction of such systems. 
We also attempt a low-rank output feedback synthesis for a linear time-varying system. 

The remainder of this paper is organized as follows. 
In the next section, we briefly summarize existing results. 
Our main contributions are presented in Section \ref{sec:main_natural_power_method}, where we provide a convergence analysis of the natural power method \eqref{eq:natural_power_method} for general square matrices. 
Section \ref{sec:ModelReduction} presents a model order reduction method based on the natural power method and an attempt at low-rank controller synthesis for a linear time-varying system. 
We conclude the paper in Section \ref{sec:conclusion} with final remarks.

\subsubsection*{Notation}
The sets of real and complex numbers are denoted by $\mathbb{R}$ and $\mathbb{C}$, respectively. The set of $n \times m$ real matrices is $\mathbb{R}^{n\times m}$. $I_n$ denotes the $n \times n$ identity matrix, and $0_{n,m}$ denotes the $n \times m$ zero matrix. For simplicity, we denote $0_{n} = 0_{n,n}$; if the dimension is trivial, $0$ is used. 
$A^{\top}$ and $A^{\dagger}$ denote the transpose and Hermitian conjugate of a matrix $A$, respectively. $\|x\|$ denotes the Euclidean norm of a vector $x$, and $\| A \| _{\rm ind}$ denotes the induced norm of a matrix $A$. For a symmetric matrix $A$, $A > 0$ ($A \geq 0$) indicates positive definiteness (semidefiniteness). $A^{1/2}$ denotes the unique positive-semidefinite square root. The eigenvalues of a square matrix $A \in \mathbb{C}^{n \times n}$ are ordered such that $|\lambda_1[A]| \geq \dots \geq | \lambda_n [A] |$. The corresponding (generalized) right eigenvector is $\psi_i [A]$, where each vector is normalized. The matrix of eigenvectors is defined as $\Psi[A] := [\psi _{1}[A] , \dots , \psi _{n}[A]]$. The Stiefel manifold is defined as $\mathrm{St}(r,n) := \{X \in \mathbb{R}^{n \times r} \mid X^{\top}X = I_r\}$.

\section{Previous Work}

Numerous principal component extraction algorithms exist for both continuous-time and discrete-time settings. 
In contrast to continuous-time frameworks, discrete-time algorithms were primarily developed to reduce computational costs. 
Notably, the NPM requires the inverse square root calculation, which may be computationally expensive. 
\cite{hua2004asymptotical} provides numerical insights into this issue. 
This paper focuses primarily on the convergence issues of the natural power method \eqref{eq:natural_power_method} rather than its computational cost, as \cite{hua1999new} noted that this method provides a standard representation for similar algorithms.

\cite{hua1999new} derived the following convergence result for symmetric positive semidefinite matrices.

\begin{prop}[{\cite[Appendix A]{hua1999new}}]\label{prop:hua1999new_convergence}
Suppose that the following two conditions hold:
\begin{enumerate}
    \item $A = A^{\top} \in \mathbb{R}^{n\times n}$ is positive semidefinite.
    \item The $r$-th eigenvalue is strictly separated from the $(r+1)$-th eigenvalue: $\lambda _{r}[A] > \lambda _{r+1}[A]$.
\end{enumerate}
If the initial condition $U[0]$ belongs to the set
\begin{align*}
    U[0] \in \Bigg{\{} 
    \Psi [A] \begin{bmatrix} K_{r} \\ K_{\perp} \end{bmatrix} 
    & 
    \in \mathrm{St}(r,n) \ \bigg| \ K_{r} \in \mathbb{R}^{r\times r}, 
    \\ &
     K_{\perp} \in \mathbb{R}^{(n-r) \times r} ,\ \mathrm{det}(K_{r}) \neq 0 \Bigg{\}},
\end{align*}
then the solution of \eqref{eq:natural_power_method} converges exponentially to the set
\begin{align}
    \mathcal{U}_{r}^{\rm (d)} := \left\{ \Psi [A]  \begin{bmatrix}
    K_{r} \\ 0_{n-r,r}
    \end{bmatrix}
    \in \mathrm{St}(r,n) \ \bigg| \ K_{r}\in \mathbb{R}^{r\times r}
    \right\} 
    \label{eq:StableEquilibriumSets_natural_power_method}
\end{align}
with a convergence rate of $\lambda _{r+1}[A]/\lambda _{r}[A]$.
\end{prop}

Since the natural power method \eqref{eq:natural_power_method} represents a standard form for this class of algorithms [\cite{hua1999new,attallah2002fast,badeau2005fast,hua2004asymptotical}], Proposition \ref{prop:hua1999new_convergence} provides an essential convergence benchmark. However, to the best of our knowledge, the convergence of this and related algorithms has not been analyzed for general square matrices.

\section{Convergence Analysis for the NPM}
\label{sec:main_natural_power_method}

In this section, we analyze the natural power method \eqref{eq:natural_power_method} for general square matrices. 
We show that the forward invariant set differs from that of the Oja flow, implying that the natural power method extracts a different subspace.

\subsection{Convergence Analysis}

This section generalizes the convergence results presented in {\cite[Appendix A]{hua1999new}} to general square matrices. The solution of the natural power method \eqref{eq:natural_power_method} remains on the Stiefel manifold $\mathrm{St}(r,n)$ provided that $U[k]^{\top} A^{\top}A U[k]$ is non-singular. This condition is always satisfied if $A$ is non-singular and $U[0] \in \mathrm{St}(r,n)$. Moreover, the following lemma holds.

\begin{lem}\label{lemma:discrete_norank_deficit}
Let $m\geq r$. If the initial condition $U[0]$ is in 
\begin{align}
    \mathcal{V}_{m,r}^{\rm (d)} := \Bigg{\{} \Psi [A] & \begin{bmatrix} K_{m,r} \\ K_{\perp} \end{bmatrix}  \in \mathrm{St}(r,n) \ \Bigg{|} \ 
    K_{m,r} \in \mathbb{C}^{m\times r},
     \nonumber \\ &
    K_{\perp} \in \mathbb{C}^{(n-m) \times r}, 
    \ \mathrm{rank}(K_{m,r}) =r \Bigg{\}}
    \label{eq:discrete_Oja_domain}
\end{align}
and the eigenvalues satisfy $| \lambda _{m}[A]| > | \lambda _{m+1}[A]|$, then $U[k] \in \mathcal{V}_{m,r}^{\rm (d)}$ and $U[k]^{\top}A^{\top}AU[k]$ is non-singular $\forall k\geq 0$.
\end{lem}

We use the simplified notation $\mathcal{V}_{r}^{\rm (d)} := \mathcal{V}_{r,r}^{\rm (d)} $.

\begin{pf}
Let $Z[k] := U[k]^{\top} A ^{\top} A U[k]$, $\Psi := \Psi [A]$, and $\Lambda := \mbox{blk-diag}( \Lambda _{m} ,\ \Lambda _{\perp} ) := \Psi ^{-1} A \Psi $, where $\Lambda _{m} \in \mathbb{C}^{m\times m} $ and $\Lambda _{\perp} \in \mathbb{C}^{(n-m) \times (n-m)}$. 
At $k=0$, the following holds:
\begin{align*}
    Z[0] =& \begin{bmatrix}
    K_{m,r} ^{\dagger} & K_{\perp}^{\dagger} 
    \end{bmatrix}
    \Lambda ^{\dagger} \Psi ^{\dagger} \Psi \Lambda 
    \begin{bmatrix}
    K_{m,r} \\ K_{\perp}
    \end{bmatrix}
    \\
    =& \begin{bmatrix}
    K_{m,r} ^{\dagger} \Lambda _{m}^{\dagger} & K_{\perp}^{\dagger} \Lambda _{\perp}^{\dagger}
    \end{bmatrix}
    \Psi ^{\dagger} \Psi 
    \begin{bmatrix}
    \Lambda _{m} K_{m,r} \\ \Lambda _{\perp}K_{\perp}
    \end{bmatrix}
    \\
    \geq &
    \lambda _{n}[ \Psi ^{\dagger} \Psi  ] K_{m,r} ^{\dagger} \Lambda _{m}^{\dagger} \Lambda _{m} K_{m,r}. 
\end{align*}
Given the assumption that $\lambda _{m}[A] \neq 0$ (i.e., $ \Lambda _{m} $ is non-singular) and $\mathrm{rank}(K_{m,r}) =r$, it follows that $Z[0] >0$, meaning $Z[0]$ is invertible. 
Therefore,
\begin{align*}
    U[1] =& A U[0] Z[0]^{-1/2}
    =
    \Psi \begin{bmatrix}
    \Lambda _{m} K_{m,r} Z[0]^{-1/2}
    \\ 
    \Lambda _{\perp}K_{\perp} Z[0]^{-1/2}
    \end{bmatrix}
\end{align*}
and $\mathrm{rank} (\Lambda _{m} K_{m,r} Z[0]^{-1/2}) =r$ since $\Lambda_m$ and $Z[0]^{-1/2}$ are invertible.

Suppose that at time $k \geq 1$, $ U[ k ] =  \Psi \begin{bmatrix}
    K_{m,r}^{(k)} \\ K_{\perp}^{(k)}
    \end{bmatrix}
    \in \mathrm{St}(r,n)$, 
where $K_{m,r}^{(k)} \in \mathbb{C}^{m\times r}$ with $\mathrm{rank}(K_{m,r}^{(k)} ) =r $ and $K_{\perp}^{(k)} \in \mathbb{C}^{(n-m) \times r}$.  
Note that $K_{\bullet } ^{(k)} = \Lambda _{\bullet} ^{k} K_{\bullet} Z[0]^{-1/2}\dots Z[ k -1] ^{-1/2}$ for $\bullet = (m,r),\perp$. 
By the same argument as for $k=0$, $Z[k]$ is invertible. Then,
\begin{align*}
    U[k +1 ] = & \Psi \begin{bmatrix}
    \Lambda _{m}^{ k +1} K_{m,r} Z[0]^{-1/2} Z[1]^{-1/2} \dots Z[ k ] ^{-1/2}
    \\ 
    \Lambda _{\perp} ^{ k +1}K_{\perp}Z[0]^{-1/2} Z[1]^{-1/2} \dots Z[ k ] ^{-1/2}
    \end{bmatrix}
    \\
    =&
    \Psi \begin{bmatrix}
    (K_{m,r}^{( k +1)}) ^{\top}
    & 
    (K_{\perp}^{( k +1)}) ^{\top}
    \end{bmatrix}^{\top}
\end{align*}
and therefore, $\mathrm{rank}(K_{m,r}^{(k +1)}) =r$. 
This implies that $Z[k+1]$ is invertible if $Z[k]$ is invertible.  
By mathematical induction, the statement holds.
\end{pf}

From Lemma \ref{lemma:discrete_norank_deficit}, we obtain the following convergence result, which parallels \cite{TsuzukiOhki2025global}.

\begin{thm}\label{thm:discrete_eigenvector_subspace}
Suppose that for a given $A\in \mathbb{R}^{n\times n}$, there exists $m \in \{r, \dots ,n-1\}$ such that $| \lambda _{m} [A] | > | \lambda _{m+1} [A] | $. 
If the initial condition $U[0] \in \mathcal{V}_{m,r}^{\rm (d)}$, where $\mathcal{V}_{m,r}^{\rm (d)}$ is defined in Eq. \eqref{eq:discrete_Oja_domain}, then the solution $U[k]$ of \eqref{eq:natural_power_method} converges to the invariant set
\begin{align}
    \mathcal{U}_{m,r}^{\rm (d)} := \Bigg{\{} \Psi [A] \begin{bmatrix} K_{m,r} \\ 0_{n-m,r} \end{bmatrix}  \in & \mathrm{St}(r,n)  \ \Bigg{|}  \ K_{m,r} \in \mathbb{C}^{m\times r},
    \nonumber \\ &
    \ \mathrm{rank}(K_{m,r}) =r \Bigg{\}}.
    \label{eq:discrete_equillibriums}
\end{align}
\end{thm}

If $m=r$, we use the simplified notation $\mathcal{U}_{r}^{\rm (d)} = \mathcal{U}_{r,r}^{\rm (d)}$. 
In this case, for any $\bar{U}_{r} \in \mathcal{U}_{r}^{\rm (d)} $, 
\begin{align}
 \bar{U}_{r}^{\top}A\bar{U}_{r} 
 =&  \begin{bmatrix} K_{r,r} \\ 0_{n-r,r} \end{bmatrix}^{\dagger} \Psi ^{\dagger} \Psi \begin{bmatrix} K_{r,r} \\ 0_{n-r,r} \end{bmatrix}  K_{r,r}^{-1} \Lambda _{r} K_{r,r} 
 \nonumber \\ 
 =& K_{r,r}^{-1} \Lambda _{r} K_{r,r},
 \label{eq:eigenvalue_preserving}
\end{align}
where $\Lambda := \mbox{blk-diag} (\Lambda _{r}, \Lambda _{\perp}) = \Psi ^{-1} A \Psi $ and $\Psi = \Psi [A]$. This means that $\bar{U}_{r}$ preserves the dominant $r$ eigenvalues of $A$ with the largest absolute values. 
Even for this case, the convergence is ensured for an invariant set, not an element of the invariant set. See Remark \ref{rem:stationary_point_NPM} below for the details. 
In practice, the spectrum gap condition in Theorem \ref{thm:discrete_eigenvector_subspace} is not easily verified for large-dimensional systems. Theorem \ref{thm:discrete_eigenvector_subspace} establishes that even if the gap condition does not hold between $r$ and $r+1$-th eigenvalues, there exists an invariant set of the NPM algorithm. 
In addition, following the argument in \cite[\S III.C]{TsuzukiOhki2025global}, the volume of $\mathcal{V}_{m,r}^{\rm (d)}$ equals that of $\mathrm{St}(r,n)$. 
Therefore, almost all initial matrices converge to $\mathcal{U}_{m,r}^{\rm (d)}$.

\begin{pf}
From the proof of Lemma \ref{lemma:discrete_norank_deficit}, $U[k] \in \mathrm{St}(r,n)$ is represented as
$U[k] = \Psi  \begin{bmatrix}
    (K_{m,r}^{(k)})^{\top} & (K_{\perp} ^{(k)})^{\top}
    \end{bmatrix} ^{\top}$, $K_{m,r}^{(k)} \in \mathbb{C}^{m\times r}$, $K_{\perp}^{(k)} \in \mathbb{C}^{(n-m) \times r}$, 
where $K_{m,r}^{(k)}$ and $K_{\perp}^{(k)}$ are defined as in the proof of Lemma \ref{lemma:discrete_norank_deficit} with $K_{m,r}^{(0)} = K_{m,r}$ and $K_{\perp}^{(0)} = K_{\perp}$. 
Let us define $V[k] := A ^{k} U[0] / | \lambda _{m}[A] |^{k} = \Psi  \begin{bmatrix}
    (\Lambda _{m}^{k} K_{m,r})^{\top} & (\Lambda _{\perp}^{k} K_{\perp} )^{\top}
    \end{bmatrix} ^{\top} / | \lambda _{m}[A] |^{k}$.  
Note that for any $k\geq 0$, $U[k] U[k]^{\top} = V[k] (V[k]^{\top} V[k])^{-1} V[k] $ holds. 
From the assumption $| \lambda _{m} [A] | >| \lambda _{m+1} [A] | $, for any $\delta \in (0, | \lambda _{m}[A]| - |\lambda _{m+1}[A]|)$, there exists $c = c(\delta ) \geq 1$ such that
\begin{align*}
    \| \Lambda _{\perp} ^{k} |\lambda _{m}[A] |^{-k} \| _{\rm ind} \leq &
    c \left| \frac{  \lambda _{m+1}[A] +\delta }{  \lambda _{m}[A]  }\right| ^{k}
    \xrightarrow{k\to \infty} 0. 
\end{align*}
This concludes that 
\begin{align*}
\lim _{k\to \infty} V[k] (V[k]^{\top} V[k])^{-1} V[k] = \Psi 
\begin{bmatrix}
* & 0_{n-m,m} \\ 0_{m,n-m} & 0 _{n-m}
\end{bmatrix}
\Psi ^{\dagger}
,
\end{align*}
and therefore, $\lim _{k\to \infty} \| K_{\perp } ^{(k)} \| _{\rm ind} =  0$. 
Since $U[k] \in \mathrm{St}(r,n)$ for any $k\geq 0$,  
$U[k]$ converges to $\mathcal{U}_{m,r}^{\rm (d)}$.
\end{pf}

\begin{figure}[!htbp]
    \centering
    \includegraphics[keepaspectratio,width=8.5cm]{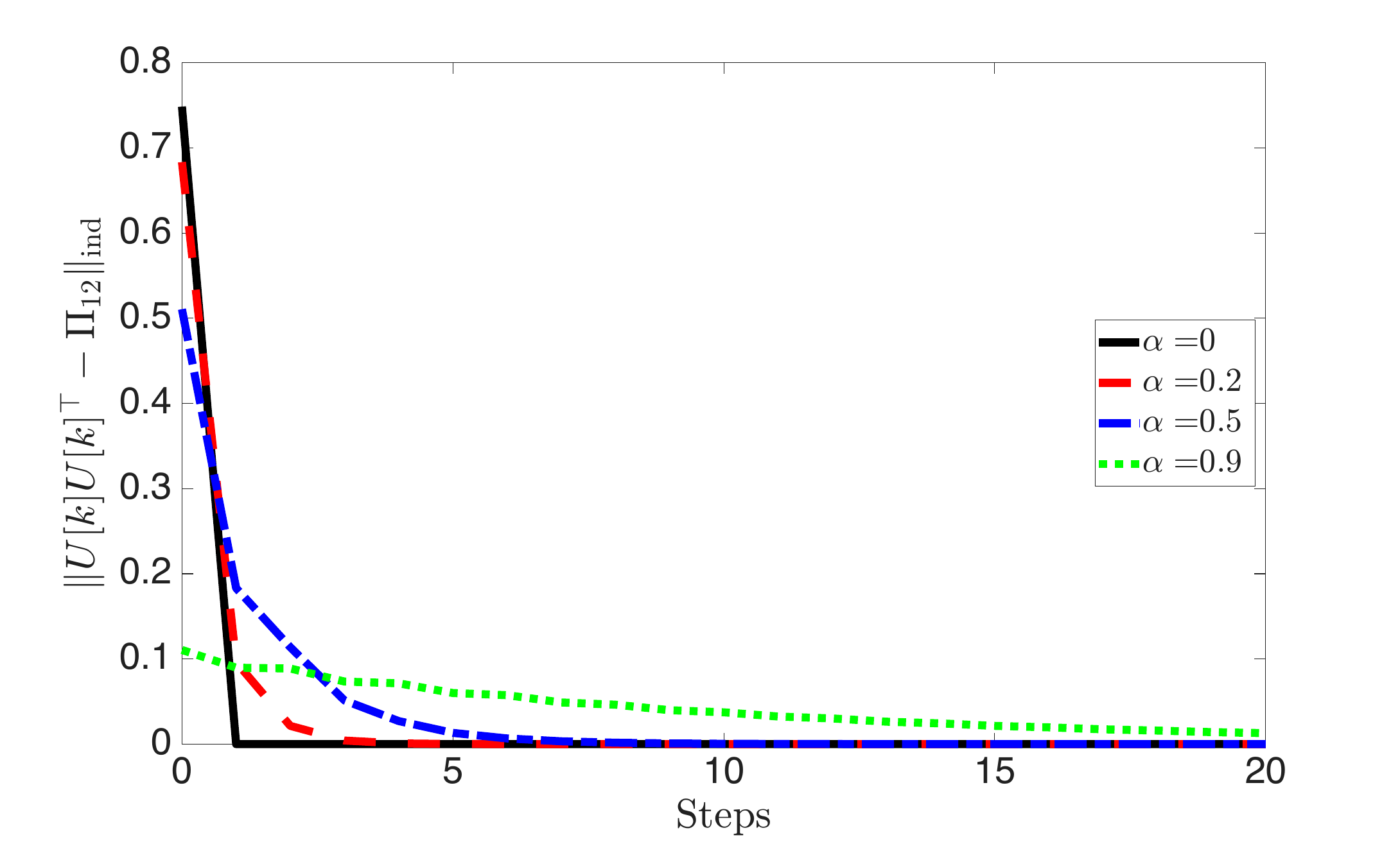}    \caption{Plot of $\| U[k] U[k] ^{\top}- \Pi _{12}\| _{\rm ind}$ for each time $k$. The black solid, red dashed, blue chain, and green dotted lines represent $\alpha =0, 0.2, 0.5, 0.9$, respectively.}
    \label{fig:example_natural_power_method_r2}
\end{figure}

From Theorem \ref{thm:discrete_eigenvector_subspace}, if the absolute values of consecutive eigenvalues, $|\lambda _{m}[A]|$ and $|\lambda _{m+1}[A]|$, are close, the convergence rate may become slow. 
Figure \ref{fig:example_natural_power_method_r2} illustrates the convergence of the natural power method \eqref{eq:natural_power_method} to $\mathcal{U}_{2}^{\rm (d)}$, where
\begin{align}
    A_{\alpha} = \begin{bmatrix}
    1 & 1 & 2
    \\
    0 & \alpha & 1
    \\
    0 & 0 & -1
    \end{bmatrix}
    ,\quad \alpha \in \{0, 0.2, 0.5, 0.9\}.
    \label{eq:A_alpha}
\end{align}
The black solid, red dashed, blue chain, and green dotted lines represent $\alpha =0, 0.2, 0.5, 0.9$, respectively.  
Note that for each $\alpha$, the eigenvectors $\psi _{i} = \psi _{i}[A_{\alpha}]$, $i=1,2,3$, are
\begin{align*}
     \psi_{1} =& \begin{bmatrix} 1 \\ 0 \\ 0 \end{bmatrix}, \quad
     \psi_{2} = \frac{1}{\sqrt{8\alpha ^{2} + 12 \alpha + 9}} \begin{bmatrix} 1+2\alpha  \\ 2 \\ -2 (1+\alpha ) \end{bmatrix} , \\
    \psi_{3} =& \frac{1}{\sqrt{2-2\alpha + \alpha^2}}\begin{bmatrix} 1, & \alpha -1 ,& 0 \end{bmatrix} ^{\top}. 
\end{align*}
The order of $\psi _{1}$ and $\psi _{2}$ is interchangeable since the absolute values of the eigenvalues are identical. 
The initial condition is set as $U[0] := X (X^{\top}X)^{-1/2}$, $X :=\begin{bmatrix} \psi _{1}+ \psi _{2} + \psi _{3}, & \psi _{2} + \psi _{3} \end{bmatrix} 
    \in \mathbb{R}^{3\times 2}$ 
for each $\alpha$.  
The distance between $\mathcal{U}_{2}^{\rm (d)}$ and $U[k]$ is defined as $\| U[k] U[k]^{\top} - \Pi _{12} \| _{\rm ind}$, where $\Pi _{12}$ is the projection onto the plane spanned by $\psi _{1}$ and $\psi _{2}$. 
Since the convergence ratio $|\lambda _{3}[A]| / |\lambda _{2}[A]| = \alpha$, convergence is slow when $\alpha =0.9$, whereas convergence is achieved in one step when $\alpha =0$.

\begin{figure}[!htbp]
    \centering
    \includegraphics[keepaspectratio,width=8.5cm]{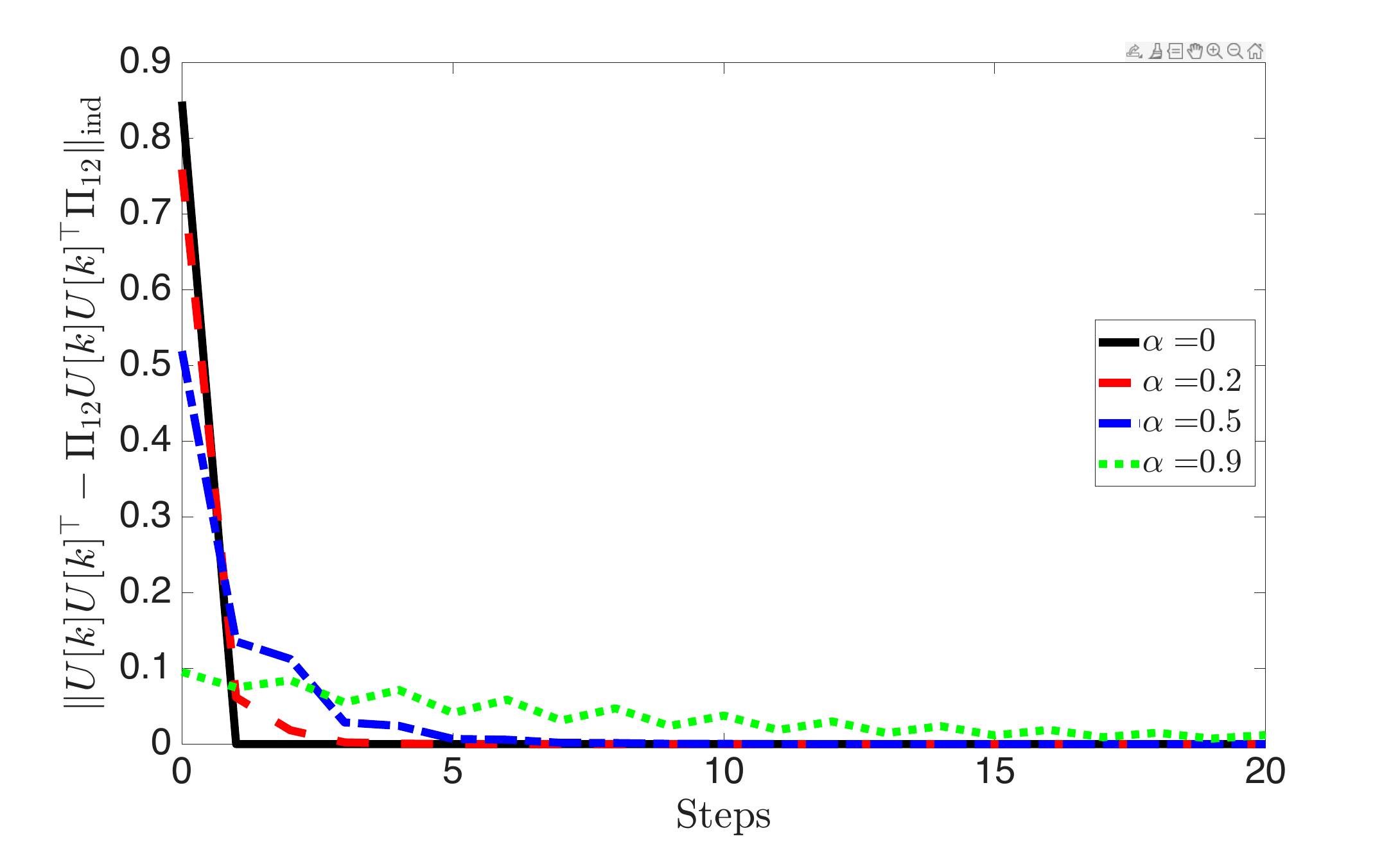}
    \caption{Plot of $\| U[k] U[k] ^{\top}- \Pi _{12} U[k] U[k] ^{\top} \Pi_{12} \|_{\rm ind}$ for each time $k$, where $U[k] \in \mathbb{R}^{3}$. The black solid, red dashed, blue chain, and green dotted lines represent $\alpha =0, 0.2, 0.5, 0.9$, respectively.}
    \label{fig:example_natural_power_method_r1}
\end{figure}

\begin{figure}[!htbp]
    \centering
    \includegraphics[keepaspectratio,width=8.5cm]{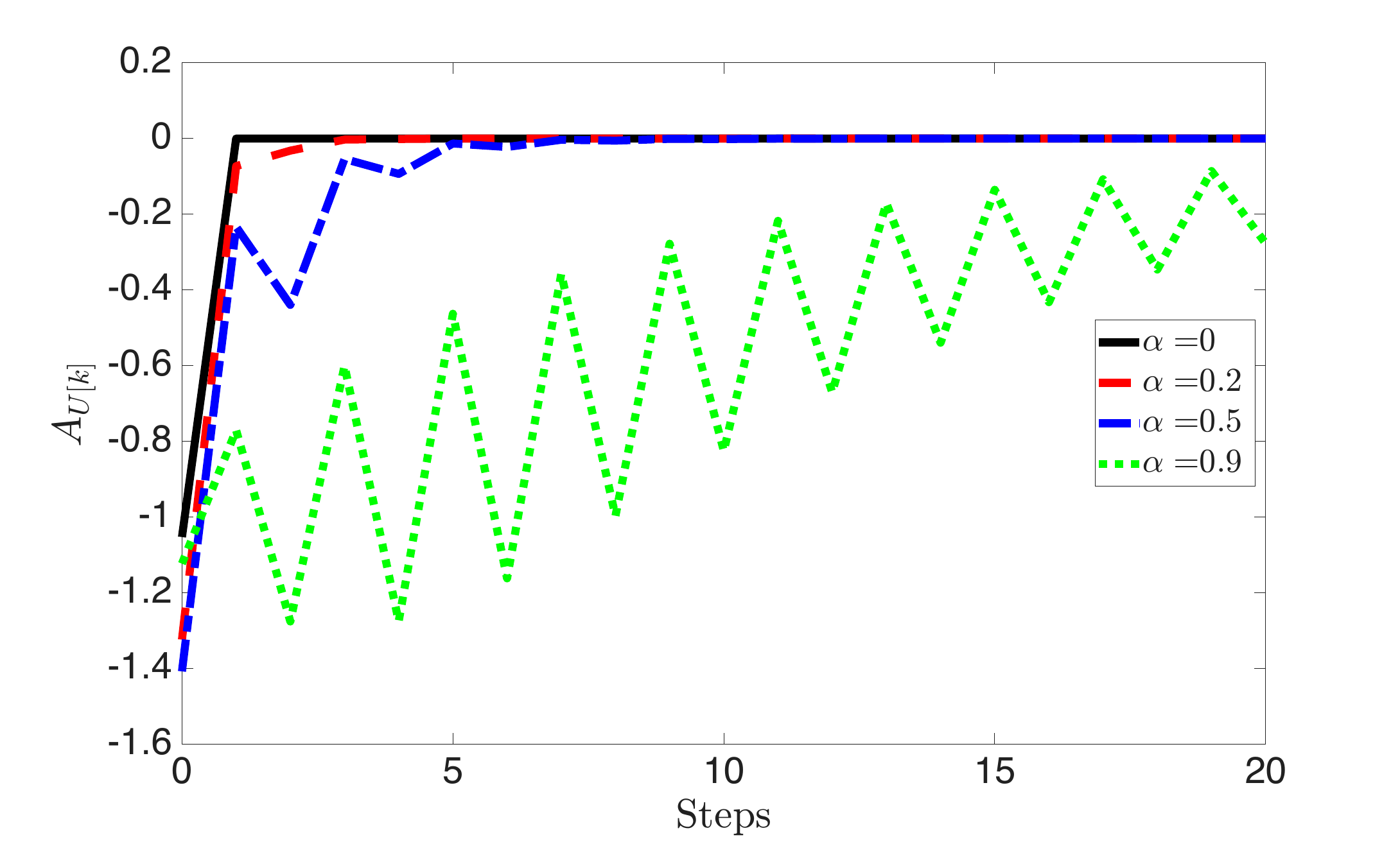}
    \caption{Plot of $A_{U[k]} = U[k]^{\top} A_{\alpha} U[k] \in \mathbb{R}$ at each time $k$. The black solid, red dashed, blue chain, and green dotted lines represent $\alpha =0, 0.2, 0.5, 0.9$, respectively.}
    \label{fig:example_natural_power_method_eigenvalues}
\end{figure}

If $r=1$, $U[k]$ converges to $\mathcal{U}_{2,1}^{\rm (d)}$ similarly to the $r=2$ case (see Figure \ref{fig:example_natural_power_method_r1}).  
However, $U[k]$ does not preserve the dominant eigenvalue of $A_{\alpha}$ (Figure \ref{fig:example_natural_power_method_eigenvalues}).  
As shown in Figure \ref{fig:example_natural_power_method_eigenvalues}, the scalar projection $A_{U[k]} \in \mathbb{R}$ converges to $0$.  
If we change the initial condition $U[0]$ using $X := \begin{bmatrix} 0.1 \psi _{1}+ \psi _{2} + \psi _{3}, & \psi _{2} + \psi _{3} \end{bmatrix}$, then $A_{U[k]}$ oscillates around $-1$ and does not converge. 

\begin{rem} \label{rem:stationary_point_NPM}
The solution $U[k]$ of \eqref{eq:natural_power_method} generally does not converge to an element of $\mathcal{U}_{m,r}^{\rm (d)}$ even if $m=r$ under the conditions of Theorem \ref{thm:discrete_eigenvector_subspace}. 
We briefly show this fact below. 
If $m=r$, any element $\bar{U}_{r}$ of $\mathcal{U}_{r}^{\rm (d)}$ can be represented by a non-singular matrix $K_{r} \in \mathbb{C}^{r\times r}$ as $\bar{U} _{r} = \Psi [A] \begin{bmatrix}
    K_{r}^{\top} & 0 _{r,n-r}
    \end{bmatrix} ^{\top} 
    =
    \Psi _{r} K_{r}$, 
where $\Psi _{r} := [\psi _{1}[A],\dots , \psi _{r} [A]]$. 
Then, 
\begin{align}
    A \bar{U}_{r} = & \Psi _{r} \Lambda _{r} K_{r}
    = 
    \Psi _{r} K_{r} K_{r}^{-1} \Lambda _{r} K_{r} =  \bar{U} _{r}A_{\bar{U}_{r}}, 
    \label{eq:exponential_property}
\end{align}
where $A_{\bar{U}_{r}} := \bar{U}_{r}^{\top}A\bar{U}_{r} $.  
If $U[k] = \bar{U}_{r} \in \mathcal{U}_{r}^{\rm (d)}$, then $U[k+1] =  \bar{U}_{r} A_{\bar{U}_{r}} (A_{\bar{U}_{r}}^{\top} A_{\bar{U}_{r}}) ^{-1/2}$.
Recall that $A_{\bar{U}_r} \in \mathbb{R}^{r\times r}$ is invertible. After $\ell$ steps, we have $U[k+ \ell ] = \bar{U}_{r} \{ A_{\bar{U}_{r}} (A_{\bar{U}_{r}}^{\top} A_{\bar{U}_{r}}) ^{-1/2} \} ^{\ell}$. 
Therefore, $\bar{U}_{r}$ is a stationary point if and only if $A_{\bar{U}_{r}} = (A_{\bar{U}_{r}}^{\top} A_{\bar{U}_{r}}) ^{1/2}$, implying that $A_{\bar{U}_{r}} =A_{\bar{U}_{r}} ^{\top} >0$, which does not hold for general matrices. 

On the other hand, the same argument as in \cite[Prop. 2]{TsuzukiOhki2024} yields that any element $\bar{U}_{r} \in \mathcal{U}_{r}^{\rm (d)}$ is an equilibrium point of the Oja flow, i.e.,
\begin{align}
    (I_{n} - \bar{U}_{r} \bar{U}_{r}^{\top} ) A \bar{U}_{r} = 0_{n,r}.
    \label{eq:equillibrium_Oja_flow}
\end{align}
This result is utilized in Section \ref{sec:ModelReduction}.
\end{rem}

Based on the arguments in Remark \ref{rem:stationary_point_NPM}, one can utilize an algorithm that possesses a stationary point.

\begin{prop}\label{prop:ensure_stationary_point}
Assume that $A$ satisfies $| \lambda _{r}[A] | > | \lambda _{r+1}[A] |$. Then, the sequence generated by the following algorithm converges to an element of $\mathcal{U}_{r}^{\rm (d)}$ if $\hat{U}[0] \in \mathcal{V}_{r} ^{\rm (d)}$ and $\hat{U}[k]^{\top}A \hat{U}[k]$ is non-singular for all $k\geq 0$:
\begin{align}
    \hat{U}[k+1] =& A \hat{U}[k] (\hat{U}[k]^{\top} A^{\top} A \hat{U}[k]) ^{-1/2} 
    \nonumber \\
    & \times (\hat{U}[k]^{\top} A^{\top} \hat{U}[k] \hat{U}[k]^{\top} A \hat{U}[k]) ^{-1/2} 
    \nonumber \\ & \times
    \hat{U}[k]^{\top} A^{\top} \hat{U}[k].
    \label{eq:ensure_stationary_point}
\end{align}
\end{prop}

\begin{pf}
Note that the difference between \eqref{eq:natural_power_method} and \eqref{eq:ensure_stationary_point} is that \eqref{eq:ensure_stationary_point} includes an additional multiplication $W[k] := (\hat{U}[k]^{\top} A^{\top} \hat{U}[k] \hat{U}[k]^{\top} A \hat{U}[k]) ^{-1/2} (\hat{U}[k]^{\top} A^{\top} \hat{U}[k]) \in \mathbb{R}^{r\times r}$, which is an orthogonal matrix.

Under the assumptions and the arguments in the proof of Lemma \ref{lemma:discrete_norank_deficit}, the solution $\hat{U}[k]$ of the algorithm \eqref{eq:ensure_stationary_point} always satisfies the non-singularity of $\hat{U}[k]^{\top} A^{\top} A \hat{U}[k]$.  
Additionally, since $W[k]$ is an orthogonal matrix, the same argument as in Theorem \ref{thm:discrete_eigenvector_subspace} implies that $U[k]$ of \eqref{eq:ensure_stationary_point} converges to $\mathcal{U}_{r}^{\rm (d)}$. 
Thus, the remainder of the proof focuses on stationarity. 

Let $\bar{U}_{r} \in \mathcal{U}_{r}^{\rm (d)}$. Then there exists a non-singular matrix $K_{r} \in \mathbb{C}^{r\times r}$ such that $\bar{U}_{r} = \Psi _{r} K_{r}$, where $\Psi _{r} := [\psi _{1}[A] ,\dots ,\psi _{r}[A]]$.  
From Eq. \eqref{eq:exponential_property}, $A \bar{U} _{r} = \bar{U}_{r} A_{\bar{U}_{r}}$, where $A_{\bar{U}_{r}} := \bar{U}_{r}^{\top}A\bar{U}_{r} \in \mathbb{R}^{r\times r}$, resulting in $\bar{U}_{r}^{\top} A^{\top} A \bar{U}_{r} = A_{\bar{U}_{r}}^{\top} A_{\bar{U}_{r}}$. 
On the other hand,
\begin{align*}
    &(\bar{U}_{r}^{\top} A^{\top} A \bar{U}_{r} ) ^{-1/2} 
    (\bar{U}_{r}^{\top} A^{\top} \bar{U}_{r} \bar{U}_{r}^{\top} A \bar{U}_{r} ) ^{-1/2} 
    \bar{U}_{r}^{\top} A \bar{U}_{r}
    \\
    =&
    (A_{\bar{U}_{r}}^{\top} A_{\bar{U}_{r}}) ^{-1} A_{\bar{U}_{r}}^{\top} = A_{\bar{U}_{r}}^{-1}
\end{align*}
holds. 
Therefore, if $\hat{U}[0] = \bar{U}_{r}  \in \mathcal{U}_{r}^{\rm (d)}$, then
    $\hat{U}[1] = A\bar{U}_{r} A_{\bar{U}_{r}}^{-1} = \bar{U}_{r} = \hat{U}[0]$.
This result concludes that the sequence generated by \eqref{eq:ensure_stationary_point} converges to an element of $\mathcal{U}_{r}^{\rm (d)}$.  
\end{pf}

Numerical experiments indicate that the assumption regarding the non-singularity of $\hat{U}[k]^{\top}A\hat{U}[k]$ is not restrictive. 
Algorithm \eqref{eq:ensure_stationary_point} requires an additional square root inverse operation, increasing computational cost. 
Since $U[k]U[k]^{\top} = \hat{U}[k] \hat{U}[k]^{\top}$, the convergence of the natural power method can be confirmed by verifying the convergence of the projection matrix $U[k]U[k]^{\top}$.

\subsection{Subspace Dimension Reduction}

We expect to efficiently increase or reduce the dimension $r$ of the extracted subspace, $\mathcal{U}_{r}^{\rm (d)}$, once an element $\bar{U}_{r} \in \mathcal{U}_{r}^{\rm (d)}$ has been computed.
Unlike the Oja flow, as discussed in \cite[\S III.D]{TsuzukiOhki2025global}, the NPM \eqref{eq:natural_power_method} is nonlinear in both $A$ and $U[k]$, which makes developing dimension-increasing algorithms challenging. We therefore focus only on reducing the dimension of the extracted subspace.

Let $r \geq 1$ and $m \in \{1,\dots ,r-1\}$ ($r>m$). Assume that the eigenvalues of $A$ are strictly separated at the $r$-th and $m$-th indices, such that $| \lambda _{r}[A] | > |\lambda _{r+1}[A] |$ and $| \lambda _{m}[A] | > |\lambda _{m+1}[A] |$.
Let $A_{\bar{U}_{r}} := \bar{U}_{r}^{\top} A \bar{U}_{r}$ be the projection of $A$ onto the $r$-dominant subspace spanned by $\bar{U}_{r}$.
We then propose the following reduced-order NPM:
\begin{align}
    u[k+1] = A_{\bar{U}_{r}} u[k] ( u[k]^{\top}A_{\bar{U}_{r}}^{\top} A_{\bar{U}_{r}} u[k]) ^{-1/2},
    \label{eq:reduced_npm}
\end{align}
where $u[0] \in \mathrm{St}(m,r)$ is the initial condition for the reduced system, and $\bar{U}_{r} u[0] \in \mathrm{St}(m,n)$ is assumed to belong to $\mathcal{V}_{m}^{{\rm (d)}}$.
From the same arguments used in Theorem \ref{thm:discrete_eigenvector_subspace}, $u[k]$ converges to the $m$-dominant subspace of the projected matrix $A_{\bar{U}_{r}}$. Consequently, the projected solution $\lim _{k\to \infty} \bar{U}_{r} u[k]$ converges to $\mathcal{U}_{m}^{\rm (d)}$.

This reduction method avoids using the original $n \times n$ matrix $A$ in subsequent iterations and is therefore computationally efficient when $r\ll n$.

\section{Applications to Control Problems}
\label{sec:ModelReduction}

This section introduces a model order reduction (MOR) and low-rank controller synthesis for discrete-time linear systems, which are the discrete-time counterparts of continuous-time ones (\cite{TsuzukiOhki2025global}). 
While related works such as \cite{hua2001optimal} considered low-rank filters for static systems and \cite{TsuzukiOhki2024a} investigated low-rank Kalman filtering for continuous-time state equations with intermittent observations, these approaches did not address the MOR of general discrete-time dynamical systems. We therefore propose low-rank methods for control problems involving discrete-time systems.

\subsection{Model Order Reduction for Discrete-Time LTI Systems}

This subsection focuses on linear time-invariant (LTI) systems. Consider the following LTI system:
\begin{align}
    x[k+1] = Ax[k] + Bu[k],\quad y[k] = Cx[k],
    \label{eq:discrete_time_LTI}
\end{align}
where $x[k] \in \mathbb{R}^{n}$ is the state, $u[k] \in \mathbb{R}^{m}$ is the control input, and $y[k] \in \mathbb{R}^{p}$ is the output. $A$, $B$, and $C$ are real matrices of appropriate dimensions. 
Throughout this subsection, we assume $|\lambda _{r}[A]| > | \lambda _{r+1}[A]|$.

Before proceeding to MOR, we present some useful preliminary results. For convenience, we define the projected matrices $A_{Y} := Y^{\top}AY$, $B_{Y} := Y^{\top}B$, and $C_{Y} := CY$ for a projection matrix $Y \in \mathbb{R}^{n\times m}$.
Equation \eqref{eq:equillibrium_Oja_flow} implies that for an orthogonal matrix $Q_{\bar{U}}:= [\bar{U}_{r},\bar{U}_{\perp}] \in \mathbb{R}^{n\times n}$, where $\bar{U}_{r} \in \mathcal{U}_{r}^{\rm (d)}$ and $\bar{U}_{\perp} \in \mathrm{St}(n-r, n)$ satisfies $(I_{n} - \bar{U}_{r} \bar{U}_{r}^{\top}) = \bar{U}_{\perp} \bar{U}_{\perp}^{\top}$, the system matrix $A$ can be simplified via similarity transformation:
\begin{align}
    Q_{\bar{U}}^{\top} A Q_{\bar{U}} = \begin{bmatrix}
    A_{\bar{U}_{r}} & \bar{U}_{r}^{\top} A \bar{U}_{\perp}
    \\
    0_{n-r,r} & A_{\bar{U}_{\perp}}
    \end{bmatrix}.
    \label{eq:similar_transformation}
\end{align}
We denote $\mathcal{U}_{r}^{\rm (d)}[A] := \mathcal{U}_{r}^{\rm (d)}$ for the invariant set of the NPM for a matrix $A$ and similarly, $\mathcal{U}_{r}^{\rm (d)}[A^{\top}]$ for the invariant set of the NPM for the matrix $A^{\top}$.  For $\bar{V}_{r} \in \mathcal{U}_{r}^{\rm (d)}[A^{\top}]$, an orthogonal matrix $Q_{\bar{V}} := [\bar{V}_{r} , \bar{V}_{\perp}] $ gives the lower block triangle matrix via similar transformation as \eqref{eq:similar_transformation}. 
Notice that from Eq. \eqref{eq:equillibrium_Oja_flow}, $A \bar{U}_{r} = \bar{U}_{r} A_{\bar{U}_{r}}$ and $\bar{V}_{r}^{\top} A = A_{\bar{V}_{r}} \bar{V}_{r}^{\top}$ hold for $\bar{U}_{r} \in \mathcal{U}_{r}^{\rm (d)}[A]$ and $\bar{V}_{r} \in \mathcal{U}_{r}^{\rm (d)}[A^{\top}]$. 
From these facts, the following property of observability and reachability Grammians holds. 
\begin{align*}
    \bar{U}_{r}^{\top} \left( \sum _{i=0}^{n-1} (A^{\top}) ^{i} C^{\top}C A^{i} \right) \bar{U}_{r} 
    =& \sum _{i=0}^{n-1} (A_{\bar{U}_{r}}^{\top}) ^{i} C_{\bar{U}_{r}}^{\top}C_{\bar{U}_{r}} A_{\bar{U}_{r}}^{i} ,
    \\
    \bar{V}_{r}^{\top} \left( \sum _{i=0}^{n-1} A ^{i} BB^{\top} (A^{\top})^{i} \right) \bar{V} _{r}
    =& \sum _{i=0}^{n-1} A_{\bar{V}_{r}} ^{i} B_{\bar{V}_{r}}B_{\bar{V}_{r}}^{\top} (A_{\bar{V}_{r}}^{\top})^{i}.
\end{align*}
These equations imply that if the system \eqref{eq:discrete_time_LTI} is observable and reachable, then the corresponding reduced system is also observable and reachable. 
However, the above observations on observability and reachability preservation seem to indicate that the MOR ensures only one of them. 
The following lemma ties the observability and reachability preservation. 
\begin{lem}\label{lem:similarity_U_V}
For any $\bar{U}_{r} \in \mathcal{U}_{r}^{\rm (d)} [A]$ and $\bar{V}_{r} \in \mathcal{U}_{r}^{\rm (d)}[A^{\top}]$, $\bar{V}_{r} ^{\top} \bar{U}_{r} $ is non-singular. Furthermore, $A_{\bar{U}_{r}} = (\bar{V}_{r} ^{\top} \bar{U}_{r} )^{-1} A_{V_{r}} \bar{V}_{r} ^{\top} \bar{U}_{r} $. 
\end{lem}
The proof is identical to \cite[Lemmas 19 and 20]{TsuzukiOhki2025global}, so we omit it. 
Lemma \ref{lem:similarity_U_V} leads the following. 
\begin{prop}\label{prop:discrete_time_mode_reduction}
Suppose that the system \eqref{eq:discrete_time_LTI} is reachable and observable.  
Let $\bar{U}_{r} \in \mathcal{U}_{r}^{\rm (d)} [A]$ and $\bar{V}_{r} \in \mathcal{U}_{r}^{\rm (d)}[A^{\top}]$. 
If the system \eqref{eq:discrete_time_LTI} is observable and reachable, then the reduced models with $(A_{\bar{U}_{r}}, B_{\bar{V}_{r}(\bar{V}_{r}^{\top}\bar{U}_{r})^{-\top}}, C_{\bar{U}_{r}})$ and $(A_{\bar{V}_{r}}, B_{\bar{V}}, C_{\bar{U}_{r}(\bar{V}_{r}^{\top}\bar{U}_{r})^{-1}})$ are observable and reachable.  
Moreover, these models share the same transfer function: 
\begin{align}
    P_{\rm rd}[z] :=& C_{\bar{U}_{r}} (zI_{r} - A_{\bar{U}_{r}}) ^{-1} B_{\bar{V}_{r}(\bar{V}_{r}^{\top}\bar{U}_{r})^{-\top}}
    \label{eq:reduced_model_transfer_function}
    \\
    =& 
    C_{\bar{U}_{r}(\bar{V}_{r}^{\top}\bar{U}_{r})^{-1}} (zI_{r} - A_{\bar{V}_{r}}) ^{-1} B_{\bar{V}_{r}}
    . \nonumber
\end{align}
\end{prop}
From Eq. \eqref{eq:eigenvalue_preserving}, if the original system \eqref{eq:discrete_time_LTI} is stable, then the reduced model \eqref{eq:reduced_model_transfer_function} is also stable. 
Therefore, the proposed MOR in Prop. \ref{prop:discrete_time_mode_reduction} preserves not only minimality but also the stability of the system. 
Further error analysis and controller synthesis based on the above MOR are left for future work.

\subsection{Applications to Singular Perturbed Systems}

Since the NPM extracts the dominant modes, it can also extract the slow modes of stable systems.  This fact leads an application to singularly perturbed systems. 
In fact, Equation \eqref{eq:similar_transformation} represents the slower modes with $A_{\bar{U}_{r}}$ and faster modes with $A_{\bar{U}_{\perp}}$. 
For the application to singularly perturbed systems, we need an algorithm that extract the fast modes, i.e., the algorithm to obtain $\bar{U}_{\perp}$ that consists of orthogonal complements of $\bar{U}_{r} \in \mathcal{U}_{r}^{\rm (d)}[A]$. 
The following lemma characterizes the properties of $\bar{U}_{\perp}$.

\begin{lem}\label{lem:fast_mode_extraction}
Suppose that $A\in \mathbb{R}^{n\times n}$ is non-singular. 
Then, the columns of $\bar{U}_{\perp} $ spans the minor left eigensubspace of $A$, i.e., $\bar{U}_{\perp} \in \mathcal{U}_{n-r}^{\rm (d)} [A^{-\top}]$. 
\end{lem}

\begin{pf}
The eigenvalues of $A^{-\top}$ are ordered as the reciprocals of the eigenvalues of $A$, such that $\lambda _{1}[A^{-\top}] = 1/\lambda _{n}[A], \dots , \lambda _{n}[A^{-\top}] = 1/\lambda _{1}[A]$. 
Because $\bar{U}_{\perp}$ satisfies $\bar{U}_{\perp}^{\top} \bar{U}_{r} = 0_{n-r,r}$, $\bar{U}_{\perp}$ has the following form: $\bar{U}_{\perp} = \Psi [A]^{-\dagger} \begin{bmatrix} 0_{n-r,r} & \tilde{K}_{\perp}^{\top} \end{bmatrix} ^{\top}$, $\tilde{K}_{\perp} \in \mathbb{C}^{(n-r) \times (n-r)}$. 
On the other hand, using $\Lambda [A] := \Psi [A]^{-1} A \Psi [A] $, 
\begin{align*}
A^{-\top}  =& \Psi [A] ^{-\dagger} \Lambda [A]^{-\dagger} \Psi [A] ^{\dagger}
= \Psi [A^{-\top}]  \Lambda [A^{-\top}] \Psi [A^{-\top}] ^{-1}
\end{align*}
holds. 
Notice that since $\Lambda [A]$ is in the Jordan form, $\Lambda [A]^{-\dagger}$ is in the same block-sized block diagonal matrix. 
Then, there exists the same block-sized block diagonal non-singular matrix $J \in \mathbb{C}^{n\times n}$ and a permutation matrix $P \in \mathbb{R}^{n\times n}$ such that $P J \Lambda [A]^{-\dagger} J^{-1}P^{\top} = \Lambda [A^{-\top}] $.  
Hence, $\Psi [A] ^{-\dagger} J^{-1} P^{\top} \Lambda [A^{-\top}] P J \Psi [A] ^{\dagger} = \Psi [A^{-\top}]  \Lambda [A^{-\top}]$ $\Psi [A^{-\top}] ^{-1}$
holds, i.e., $\Psi [A^{-\top}] ^{-1} P J = \Psi [A] ^{-\dagger} $. Therefore, 
\begin{align*}
\bar{U}_{\perp} =& \Psi [A^{-\top}] ^{-1} P J \begin{bmatrix} 0_{r,n-r} \\ \tilde{K}_{\perp} \end{bmatrix}
=
\Psi [A^{-\top}] ^{-1} \begin{bmatrix} K_{\perp} \\ 0_{r,n-r} \end{bmatrix}
\end{align*}
holds for some $K_{\perp} \in \mathbb{C}^{(n-r)\times (n-r)}$, i.e., $\bar{U}_{\perp} \in \mathcal{U}_{n-r}^{\rm (d)}[A^{-\top}]$. 
\end{pf}

From Lemma \ref{lem:fast_mode_extraction}, the NPM for $A^{-\top}$ provides the fast modes extraction by $\bar{U}_{\perp}  \in \mathcal{U}_{n-r}^{\rm (d)} [A^{-\top}] $. 
From \eqref{eq:similar_transformation}, the state state equation \eqref{eq:discrete_time_LTI} is decomposed as 
\begin{align*}
x_{s} [k+1] =& A_{\bar{U}_{r}} x_{s}[k] + \bar{U}_{r}^{\top} A \bar{U}_{\perp} x_{f}[k] + \bar{U}^{\top}B u[k]
\\
x_{f} [k+1] = & A_{\bar{U}_{\perp}} x_{f}[k] + \bar{U}_{\perp}^{\top}B u[k]
\end{align*}
with the slow mode $x_{s} :=  \bar{U}_{r}^{\top} x$ and the fast mode $x_{f} :=  \bar{U}_{\perp}^{\top} x$, which is the standard form of singularly perturbed systems. 
It is worth noting that since NPM is a recursive algorithm, this NPM-based extraction of slow and fast modes is applicable to linear slowly time-varying systems.

\subsection{Reduced Rank Controller Synthesis for LTV Systems}

In this subsection, we present a low-rank observer-based controller synthesis to linear time-varying (LTV) systems. 
For continuous-time systems, \cite{tranninger2022detectability} demonstrated low-rank filtering for a class of LTV systems under the infinitesimal convergence time condition of a variation of the Oja flow. 
As a first step, consider the following slowly varying LTV system: 
\begin{align}
    x[k+1] = A[k] x[k] + B u[k],\quad 
    y[k] = C x[k], 
    \label{eq:LTV}
\end{align}
where $A[k] := Q[k]^{\top} A_{\alpha} Q[k]$, $B := \begin{bmatrix} 0 & 0 & 1 \end{bmatrix}^{\top}$
, $C := \begin{bmatrix} 1 & 0 & 0 \end{bmatrix}$,
\begin{align*}
    Q [k] :=& Q_{1}[k]Q_{2}[k], \ 
    Q_{1}[k] :=
    \begin{bmatrix} \cos (\omega_{1}k) & -\sin (\omega_{1}k) & 0
    \\
    \sin (\omega_{1}k) & \cos(\omega_{1}k) & 0 \\ 0& 0 & 1 
    \end{bmatrix},
    \\
    Q_{2}[k]:= &
    \begin{bmatrix}1 & 0 & 0 
    \\ 0 & \cos (\omega_{2}k) & -\sin (\omega_{2}k) 
    \\
    0 & \sin (\omega_{2}k) & \cos(\omega_{2}k)  \end{bmatrix},
\end{align*}
with $\omega _{1}=0.01$, $\omega_{2}=0.03$, and $A_{\alpha}$ defined as in \eqref{eq:A_alpha}. 

To stabilize the system \eqref{eq:LTV}, we consider the following observer-based controller:
\begin{align}
    \hat{x}[k+1] =& A[k] \hat{x}[k] + B u[k] + U[k] L[k] (y[k] - C \hat{x}[k]),
    \\
    u[k] =& F[k] V[k]^{\top} \hat{x}[k],
\end{align}
where $U[k] \in \mathrm{St}(2,3)$ and $V[k] \in \mathrm{St}(2,3)$ are the solutions of the NPM \eqref{eq:natural_power_method} for $A[k]$ and $A[k]^{\top}$, respectively.  
$L[k] \in \mathbb{R}^{2 \times 1}$ and $F[k] \in \mathbb{R}^{1\times 2}$ are the observer and feedback gains, respectively. 

\begin{figure}[!htbp]
    		\centering
    		\includegraphics[keepaspectratio,width=8.5cm]{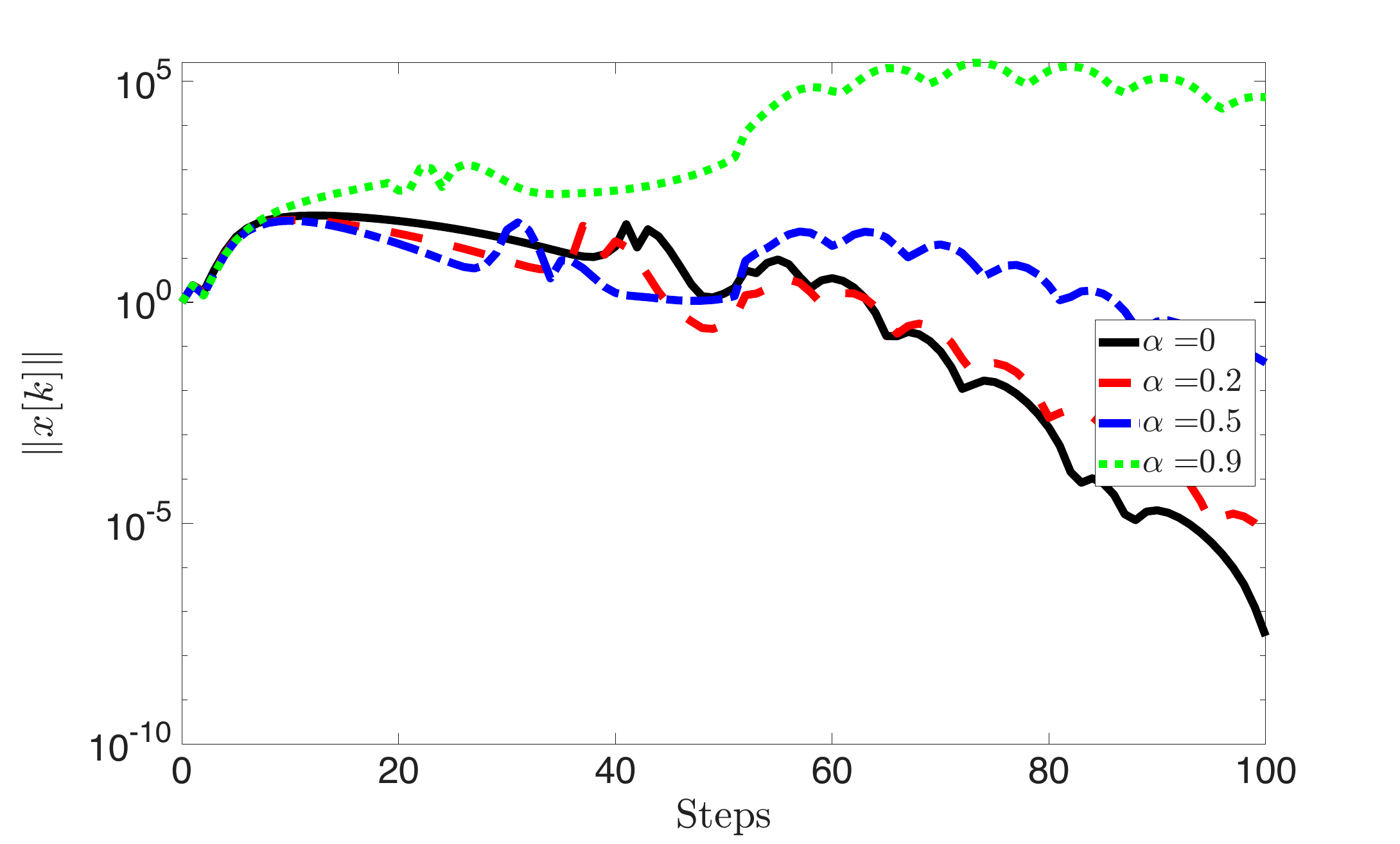}
    		\caption{Plot of $\| x[k] \| $ for each time $k$. The black solid, the red dashed, the blue chain, and green dotted lines represent $\alpha =0,0.2,0.5,0.9$, respectively.  
		}
		\label{fig:example_natural_power_method_LTV_control}
		\end{figure}

\begin{figure}[!htbp]
    		\centering
    		\includegraphics[keepaspectratio,width=8.5cm]{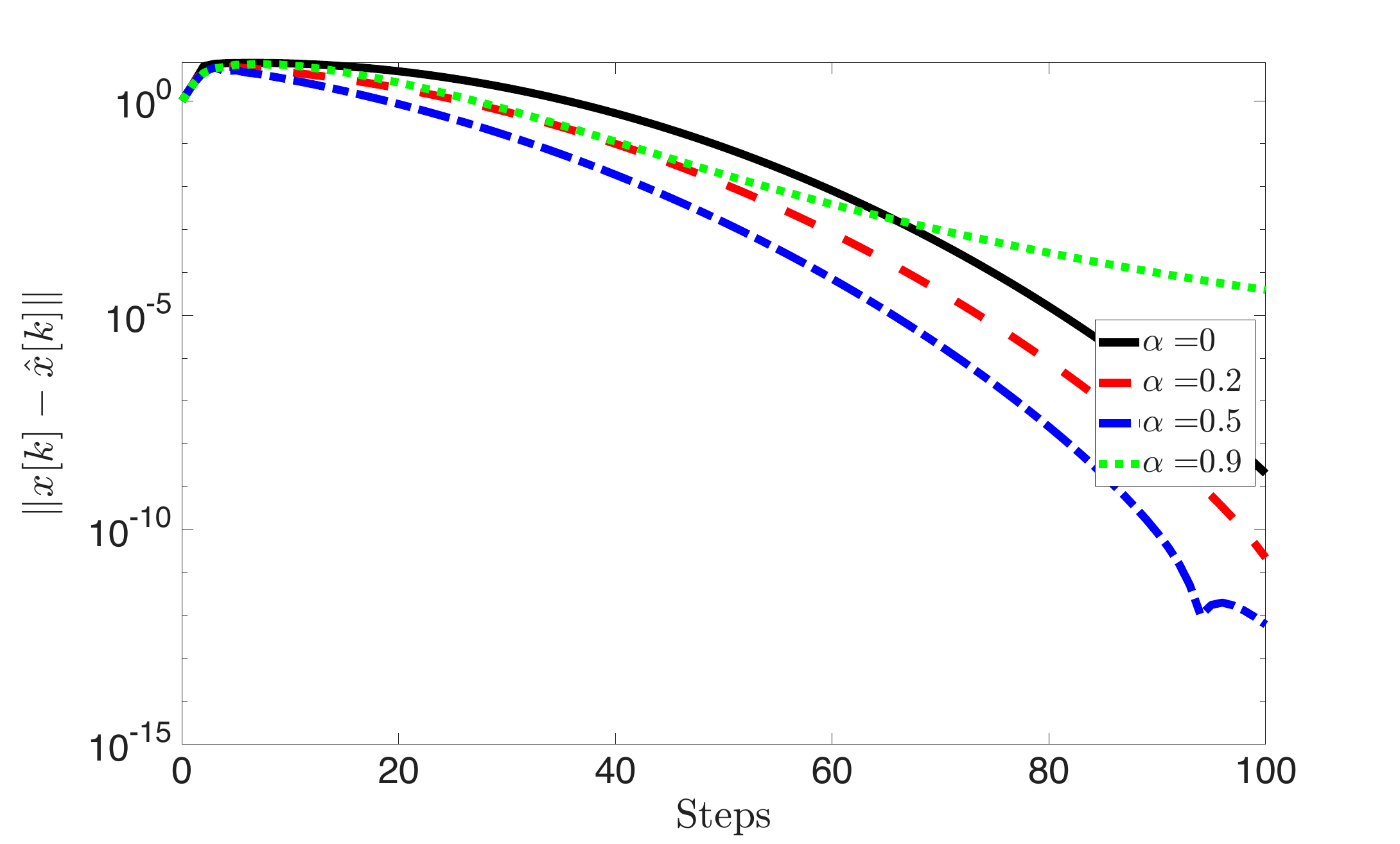}
    		\caption{Plot of $\| x[k] - \hat{x}[k]\| $ for each time $k$. The black solid, the red dashed, the blue chain, and green dotted lines represent $\alpha =0,0.2,0.5,0.9$, respectively.  
		}
		\label{fig:example_natural_power_method_LTV_estimation}
		\end{figure}

To the best of the authors' knowledge, there is no comprehensive synthesis of feedback and observer gains for LTV systems. 
In this paper, we apply the pole assignment method to obtain the gains instantaneously at each step. Specifically, at each step $k$, $L[k]$ and $F[k]$ are designed such that $A_{U[k]} - L[k] C_{U[k]} \in \mathbb{R}^{2\times 2}$ and $A_{V[k]} - B_{V[k]} F[k] \in \mathbb{R}^{2\times 2}$ have eigenvalues at $0.5$ and $0.7$. 
Figures \ref{fig:example_natural_power_method_LTV_control} and \ref{fig:example_natural_power_method_LTV_estimation} show the trajectory of $\| x[k]\|$ and the estimation error of the observer for each $\alpha$, respectively. 
The figures indicate that state estimation and feedback control function effectively for $\alpha =0, 0.2, 0.5$.  
However, the controller does not perform well for $\alpha = 0.9$, where the convergence rate of the NPM is slow. 
In this case, the proposed method does not propose a suitable feedback controller. 
Further analysis, particularly on the tracking capabilities of the NPM for slowly-varying LTV systems, is required.

\section{Conclusion}
\label{sec:conclusion}

	We investigated the convergence of the discrete-time natural power method, demonstrating its ability to extract the dominant subspace corresponding to the $r$ eigenvalues with the largest absolute values—a crucial distinction from the Oja flow. We applied this result to develop methods for MOR and low-rank controller synthesis for discrete-time systems. The MOR scheme was shown to preserve key system properties. Future work will focus on the rigorous analysis of error bounds and broader applications to time-varying control and estimation problems.

\end{document}